\documentclass{amsart}

\usepackage{t1enc}
\usepackage{fancyhdr}  

\usepackage{graphicx}
\usepackage{hyperref}

\hypersetup{
  colorlinks   = true,   
  urlcolor     = blue,   
  linkcolor    = blue,   
  citecolor   = red      
}

\usepackage[citation-order, msc-links]{amsrefs}

\issueinfo{10}{1}{number 1}{2016}
\copyrightinfo{2016}{\href{http://albanian-j-math.com}{Albanian Journal of Mathematics}}
\def\issn{{\sc \textbf{ISSN: }} 1930-1235; }
\def\issueyear{\textbf{2016}}
\PII{\issn  (\issueyear)}

\pagespan{47}{80}

\begin{document}

\title{Application of Orthonormal Bernoulli Polynomials for Approximate Solution of Some Volterra Integral Equations}


\author{Udaya Pratap Singh}
\address{Department of Applied Sciences, Rajkiya Engineering College, Sonbhadra \\  Uttar Pradesh, India}
\curraddr{}
\email{upsingh1980@gmail.com}
\thanks{}


\dedicatory{Rajkiya Engineering College, Sonbhadra, U. P., India}

 


{\let\thefootnote\relax\footnotetext{\textit{English title:} \textbf{Application of Orthonormal Bernoulli Polynomials for Approximate Solution of Some Volterra Integral Equations}.}}

\maketitle


\hrule

\vspace{.2cm}

\begin{abstract}
In this work, a new approach has been developed to obtain numerical solution of linear Volterra type integral equations by obtaining asymptotic approximation to solutions. Using the classical Bernoulli polynomials, a set of orthonormal polynomials have been derived, and these orthonormal polynomials have been used to form an operational matrix of integration which is has been implemented to find numerical or exact solution of non-singular Volterra integral equations. Two linear Volterra integral and two convolution integral equations of second kind have been solved to demonstrate the effectiveness of present method. Obtained approximate solutions have been compared with the exact solutions for numerical values. High degree of accuracy of numerical solutions has established the credibility of the present method.

\end{abstract}

\vspace{.4cm}

\hrule

\vspace{.4cm}

\noindent \text{Mathematics Subject Classes 2010:} 45A05; 45D05; 45L05; 65R20  

\noindent \text{Keywords:} Volterra integral equation; Bernoulli polynomials; orthonormal polynomials

\vspace{.4cm}

\hrule

\vspace{.4cm}



\section{Introduction}

Many physical problems are formulated as integral equations. Diffusion problems, heat conduction, concrete problem of physics and mechanics, unsteady Poiseuille flow in a pipe are some such examples. Also, such integral equations arise natural way in different applications of potential theory, continuum mechanics, electricity and magnetism, geophysics, antenna, synthesis problem, population genetics communication theory, mathematical modelling of economics, radiation problems, fluid mechanics, problems of astrophysics concerning transport of particles, and many more. Bulk of literature is available on Volterra and Fredholm integral equations \cite{Atkinson2011,Yousefi2006,Maleknejad2007,Samadyar2019,Bhrawy2012}. Bernoulli polynomials and its properties have been also discussed by many authors \cite{Cheon2003,Kurt2011,Natalini2003}.

Volterra integral equations uncover several difficulties referring to mathematical physics such as heat conduction difficulties. In recent years, researchers have focused their attention to find approximate solutions of integral equations. Xu \cite{Xu2007} adopted method of variational iteration, Pandey, et. al. \cite{Pandey2009} applied homotopic perturbation and method of collocation. Cheon \cite{Cheon2003} discussed possible applications of Bernoulli polynomials and functions in numerical analysis. Some other latest investigations include uses of Chebyshev polynomials \cite{Maleknejad2007}, Legendre polynomials \cite{Nemati2015}, Laguerre polynomials and Wavelet Galerkin method \cite{Rahman2012}, Legendre wavelets \cite{Yousefi2006}, the operational matrix \cite{Sahu2019}, Bernoulli matrix method \cite{Tohidi2013}. Recently, Bernoulli polynomials were used by Tohidi and Khorsand \cite{Bhrawy2012,Tohidi2013a} to solve second-order linear system of partial differential equations, Mohsenyzadeh \cite{Mohsenyzadeh2016} to solve linear Volterra integral equations, and Samadyar and Mirazee \cite{Samadyar2019} to find numerical solution for singular partial integro-differential equation of fractional order.

However, the numerical methods have certain limits and, therefore, there is always a need for an efficient method to produce more accurate numerical solution of integral equations.

In this work, it is proposed to introduce a new operational matrix of integration for orthonormal polynomials to reduce Volterra type integral equations into a system of algebraic equations. The operational matrix and method is a refinement of that used by Singh et al. \cite{Singh2019}. By using operational matrix of these orthonormal polynomials, exact solution for many Volterra integral equations can be obtained. Furthermore, the solutions to the integral equations solved with present method have been compared with exact solution of the problem. 

\section{Bernoulli Polynomials}
The monic polynomials
\begin{equation}\label{eq.1 : Basic Bernoulli Polynomials}
{B_n}(\zeta ) = \sum\limits_{j = 0}^n {\,\left( {\begin{array}{*{20}{c}}
		n \\ 
		j 
		\end{array}} \right)\,\,{B_j}(0)\,{\zeta ^{n - j}},\,\,\,\,\,\,\,n = 0,1,2,...} \,\,;\,\,\,\,\,0 \le \zeta  \le 1
\end{equation}
were introduced by Jacob Bernoulli in early sixteenth century, where $B_k(0)$  are the Bernoulli numbers. To have a better understanding, first few Bernoulli polynomials are represented as:
\begin{equation} \label{eq.2 : Bernoulli 1}
{B_0}(\zeta ) = 1
\end{equation}
\begin{equation} \label{eq.3 : Bernoulli 2}
 B_1(\zeta)=\zeta-\frac{1}{2}
\end{equation}
\begin{equation} \label{eq.4 : Bernoulli 3}
B_2(\zeta)=\zeta^2-\zeta+\frac{1}{6}
\end{equation}
\begin{equation} \label{eq.5 : Bernoulli 4}
B_3(\zeta)=\zeta^3-\frac{3}{2}\zeta^2+\frac{1}{2}\zeta
\end{equation}
\begin{equation} \label{eq.6 : Bernoulli 5}
B_4(\zeta)=\zeta^4-2\zeta^3+\zeta^2-\frac{1}{30}
\end{equation}
\begin{equation} \label{eq.7 : Bernoulli 6}
{B_5}(\zeta ) = {\zeta ^5} - {\textstyle{5 \over 2}}{\zeta ^4} + {\textstyle{5 \over 3}}{\zeta ^3} - {\textstyle{1 \over 6}}\zeta 
\end{equation}
\begin{equation} \label{eq.8 : Bernoulli 7}
{B_6}(\zeta ) = {\zeta ^6} - 3{\zeta ^5} + \frac{5}{2}{\zeta ^4} - \frac{1}{3}{\zeta ^2} + \frac{1}{{42}}
\end{equation}

However, the name \textit{Bernoulli Polynomials} was coined by J. L. Raabe in 1851, a thorough study of these polynomials for arbitrary value of its variable   was first done by Leonhard Euler in 1755, who showed in his book "\textit{Foundations of differential calculus}" that these polynomials satisfy the finite difference relation.

\begin{equation}\label{eq.9 : Bernoulli recurrence formula}
B_n(\zeta+1)-B_n(\zeta)=n\zeta^{n-1},n\ge1
\end{equation}
Bernoulli Polynomials form a complete basis over $[0,1]$ \cite{Kreyszig1978} and can also be extracted from its generating function  
\begin{equation}\label{eq.10 : GF of Bernoulli Polynomials}
\frac{\gamma e^{\zeta\gamma}}{e^\gamma-1}=\sum_{n=0}^{\infty}{B_n(\zeta)\frac{\gamma^n}{n!}}\left(\left|\zeta\right|<2\pi\right)
\end{equation}
Some interesting properties of Bernoulli polynomials \cite{Costabile2001} are :
\begin{equation}\label{eq.11 : properties of Bernoulli Polnmls}
\left. \begin{gathered}
B'_n(\zeta ) = n{B_{n - 1}}(\zeta ),\,\,\,n \ge 1 \hfill \\
\int_0^1 {{B_n}(z)dz = 0,\,\,\,\,\,\,\,\,\,\,n \ge 1}  \hfill \\
{B_n}(\zeta  + 1) - {B_n}(\zeta ) = n{\zeta ^{n - 1}},\,\,n \ge 1 \hfill \\ 
\end{gathered}  \right\}
\end{equation}

For more properties and generalizations such as derivative, integration and differential equation of Bernoulli polynomials can be found in the significant works \cite{Kurt2011,Natalini2003,Lu2011,Costabile2006,Momenzadeh2017}.

\section{The Orthonormal Polynomials}
It can be easily verified that the polynomials ${B_n}(x)\,(n \ge 1)$  given by eq. (\ref{eq.1 : Basic Bernoulli Polynomials})  are orthogonal to $B_o(x)$ with respect to standard inner product on  ${L^2} \in [0,1]$. Using this property, an orthonormal set of polynomials can be derived for any $B_n$ with Gram-Schmidt orthogonalization. First ten orthonormal polynomials derived for $B_9(x)$ : 

\begin{equation} \label{eq. 12 : Ortho 1}
{\phi _{0\,}}(\zeta ) = 1
\end{equation}
\begin{equation} \label{eq. 13 : Ortho 2}
{\phi _1}(\zeta ) = \sqrt 3 ( - 1 + 2\zeta )
\end{equation}
\begin{equation} \label{eq. 14 : Ortho 3}
\phi_2\left(x\right)=\sqrt5\left(1-6x+6x^2\right)
\end{equation}
\begin{equation} \label{eq. 15 : Ortho 4}
\phi_3(\zeta)=\sqrt7(-1+12\zeta-30\zeta^2+20\zeta^3)
\end{equation}
\begin{equation} \label{eq. 16 : Ortho 5}
\phi_4(\zeta)=3(1-20\zeta+90\zeta^2-140\zeta^3+70\zeta^4)
\end{equation}
\begin{equation} \label{eq. 17 : Ortho 6}
\phi_5(\zeta)=\sqrt{11}(-1+30\zeta-210\zeta^2+560\zeta^3-630\zeta^4+252{\zeta5}^)
\end{equation}
\begin{equation} \label{eq. 18 : Ortho 7}
{\phi _6}(\zeta ) = \sqrt {13} \left( \begin{array}{l}
1 - 42\zeta  + 420{\zeta ^2} - 1680{\zeta ^3}\\
+ 3150{\zeta ^4} - 2772\zeta {}^5 + 924{\zeta ^6}
\end{array} \right)
\end{equation}
\begin{equation} \label{eq. 19 : Ortho 8}
{\phi _7}\left( x \right) = \sqrt {15} \left( \begin{array}{l}
- 1 + 56x - 756{x^2} + 4200{x^3}\\
- 11550{x^4} + 16632{x^5} - 12012{x^6} + 3432{x^7}
\end{array} \right)
\end{equation}
\begin{equation} \label{eq. 20 : Ortho 9}
{\phi _8}\left( x \right) = \sqrt {17} \left( \begin{array}{l}
- 1 + 72x - 1260{x^2} + 9240{x^3} - 34650{x^4}\\
+ 72072{x^5} - 84084{x^6} + 51480{x^7} - 12870{x^8}
\end{array} \right)
\end{equation}
\begin{equation} \label{eq. 21 : Ortho 10}
{\phi _9}\left( x \right) = \sqrt {19} \left( \begin{array}{l}
- 1 + 90x - 1980{x^2} + 18480{x^3} - 90090{x^4} + 252252{x^5}\\ - 420420{x^6} + 411840{x^7} - 218790{x^8} + 48620{x^9}
\end{array} \right)
\end{equation}

\section{Approximation of Functions}
Let $\phi=\left\{\phi_0,\phi_1,\phi_2,...,\phi_n\right\}$ contains first n+1 orthonormal polynomials derived for Bernoulli polynomial $B_n(x)$. Because $\phi\subset L^2[0,1]$ and $span\{\phi\}$ is a finite dimensional space, any function $f\in L^2[0,1]$ has a unique and best approximation $\hat{f}\in span\{\phi\}$ such that $\forall g\in span\{\phi\}, ||f-f||\le ||f-g||$, and 
\begin{equation}\label{eq.22 : series approximation}
f = \hat f = \mathop {\lim }\limits_{n \to \infty } \sum\limits_{k = 0}^n {\,{c_k}\,{\phi _k}(\zeta )}
\end{equation}
where $c_k=\left\langle f\middle|\phi_k\right\rangle$, and $\left\langle.\middle|.\right\rangle$ is the standard inner product on $L^2\in[0,1]$ \cite{Mahmoud2018}.

For numerical approximation, series in eq.(\ref{eq.22 : series approximation}) can be truncated after certain number of terms $-$ say $n=m$ terms, so that:
\begin{equation}\label{eq.23 : finite series approximation}
f(\zeta ) \cong \,\,\sum\limits_{k = 0}^m {\,{c_k}\,{\phi _k} = {C^T}\,\phi (\zeta )}
\end{equation}

where $C=\left(c_0,c_1,c_2,...,c_m\right), \phi(\zeta)=\left(\phi_0,\phi_1,\phi_2,...,\phi_m\right)$ are column vectors, and number of terms $m$ is chosen to meet required accuracy.

\section{Construction of operational matrix}
The orthonormal polynomials (as shown in eq. (\ref{eq. 12 : Ortho 1}-\ref{eq. 21 : Ortho 10})) can be expressed as:
\begin{equation}\label{eq.24 : Integration of Ortho Pol. of deg. 0}
\int_0^\zeta  {{\phi _o}(\eta )d\eta }  = {\phi _o}(\zeta ) + \frac{1}{{2\sqrt 3 }}{\phi _1}(\zeta )
\end{equation}

\begin{equation}\label{eq.25 : Integration of Ortho Pols.}
\begin{array}{l}
\int\limits_0^\zeta  {{\phi _i}(x)dx = } \,\,\,\,\,\frac{1}{{2\sqrt {(2i - 1)(2i + 1)} }}{\phi _{i - 1}}(\zeta )\\
\,\,\,\,\,\,\,\,\,\,\,\,\,\,\,\,\,\,\,\,\,\,\,\,\,\, + \frac{1}{{2\sqrt {(2i + 1)(2i + 3)} }}{\phi _{i + 1}}(\zeta ),\,\,\,\,(\,for{\rm{ }}\,i = 1\,,2,...\,,m)
\end{array}
\end{equation}

Relations (\ref{eq.24 : Integration of Ortho Pol. of deg. 0}-\ref{eq.25 : Integration of Ortho Pols.}) can be represented in combined form as:
\begin{equation}\label{eq.26 : combined expression for Ortho Pols. Integrated}
\int\limits_0^\zeta  {{\phi}(\eta )d\eta  = \,\,} {\Theta _{(m + 1)}}\,{\phi}(\zeta )
\end{equation}

where $\zeta\in[0,1]$ and $\Theta_{m+1}$ is operational matrix of order $(m+1)$ given as :
\begin{equation} \label{eq.27 : operational matrix}
{\Theta _{m + 1}}\, = \frac{1}{2}\left[ {\begin{array}{*{20}{c}}
	1&{\frac{1}{{\sqrt {1.3} }}}&0& \cdots &0\\
	{\frac{{ - 1}}{{\sqrt {1.3} }}}&0&{\frac{1}{{\sqrt {3.5} }}}& \cdots &0\\
	0&{\frac{{ - 1}}{{\sqrt {3.5} }}}&0& \ddots & \vdots \\
	\vdots & \vdots & \ddots &0&{\frac{1}{{\sqrt {(2m - 1).(2m + 1)} }}}\\
	0&0& \cdots &{\frac{{ - 1}}{{\sqrt {(2m - 1).(2m + 1)} }}}&0
	\end{array}} \right]
\end{equation}

\section{Solution of Linear Volterra Integral Equations }
Consider the linear Volterra integral equation of second kind:
\begin{equation} \label{eq. 28 : Volterra Equation}
y(\zeta ) = f(\zeta ) + \int\limits_0^\zeta  {\kappa (\zeta ,x)\,\,y(x)dx,\,\,\,\,0 \le \zeta  \le 1}
\end{equation}

where $y(\zeta)$ is some real valued function, $f(\zeta)$ and $k(\zeta,x)$ are continuous functions defined on $I=[0,1]$ and $S = \left\{ {\left( {\zeta ,\,\,x} \right)\,:0 \le x \le \zeta  \le 1} \right\}$  respectively. Following the classical theory of Volterra integral equations, eq. (\ref{eq. 28 : Volterra Equation}) possesses a unique solution in $C\left[ {0,1} \right]$. Moreover, if  $f(\zeta)$ and $k(\zeta,x)$ are continuously $n-differential$ on $[0, 1]$ and $S$ respectively, the unique solution of eq. (\ref{eq. 28 : Volterra Equation}) is also continuously  $n-differential$ on $[0, 1]$.

Representing $y(\zeta)$  and $f(\zeta)$  as :
\begin{equation} \label{eq. 29 : Transform of y}
y(\zeta ) = {C^T}\,\phi (\zeta )
\end{equation}
\begin{equation} \label{eq. 30 : Transform of f}
f(\zeta ) = {F^T}\,\phi (\zeta ),
\end{equation}

eq. (\ref{eq. 28 : Volterra Equation}) can be re-written as:
\begin{equation} \label{eq. 31 : Transform of full eq.}
{C^T}\phi (\zeta ) = {F^T}\phi (\zeta ) + {C^T}\int\limits_0^\zeta  {\kappa (\zeta ,x)\,\,\phi (x)dx}  = {F^T}\phi (\zeta ) + {C^T}\,{\Phi _{m + 1}}\,\,\phi (\zeta )
\end{equation}
which gives
\begin{equation} \label{eq. 32 : C in terms of Operational Matrix}
{C^T} = {\left( {I - {\Phi_{m + 1}}} \right)^{ - 1}}{F^T}
\end{equation}

where, ${\Phi _{m + 1}}\,\,\phi (\zeta ) = \int\limits_0^\zeta  {\kappa (\zeta ,x)\,\,\phi (x)dx}$, and ${\Phi _{m + 1}}$  is associate matrix of ${\Theta _{m + 1}}$ of order $m+1$, for illustration, it can be readily observed from eq. (\ref{eq.26 : combined expression for Ortho Pols. Integrated}), that  ${\Phi_{m + 1}} = c\,{\Theta _{m + 1}}$ if $\kappa (\zeta ,x) = c$ (constant) and ${\Phi _{m + 1}} = \,\Theta _{m + 1}^j$ if $\kappa (\zeta ,x) = {(\zeta  - x)^j}, (j>0)$.

\section{Error Estimate and Convergence Analysis}
\subsection*{Theorem 1}
Suppose $y(\zeta)$  be a defined and continuous for  $\zeta  \in [0,1]$,  $\phi _n^y(\zeta ) = \sum\limits_{n = 0}^\infty  {{c_k}{\phi _k}} $ be an approximation of  $y(\zeta)$  in terms of orthonormal Bernoulli polynomials $({\phi _k})$ , and $R_n(\zeta)$  be the remainder due truncation, then following relations hold.

\begin{equation} \label{eq. 33}
\phi _n^y(\zeta ) = y(\zeta ) + {R_n}(\zeta );\,\,\,\forall \,\,x \in [0,1]
\end{equation}
\begin{equation} \label{eq. 34}
\phi _n^y(\zeta ) = \int\limits_0^1 {y(\eta )d\eta  + \sum\limits_{k = 0}^n {\frac{{{\phi _k}(\zeta )}}{k}\left( {{y^{(k - 1)}}(1) - {y^{(k - 1)}}(0)} \right)} } 
\end{equation}
\begin{equation} \label{eq. 35}
{R_n}(\zeta ) =  - \frac{1}{{n!}}\int\limits_0^1 {\phi _n^ * (\zeta  - \eta ){y^{(n)}}(\eta )d\eta } 
\end{equation}
where $\phi _n^ * (\zeta ) = {\phi _n}(\zeta  - [\zeta ])$  and $\left[ \cdot \right] $  is the greatest integer function.
\subsection*{Proof} See Tohidi and Kiliçman \cite{Tohidi2013} or Mahmoud \cite{Mahmoud2018}.

\subsection*{Theorem 2}
Suppose that $y(\zeta ) \in {C^\infty }[0,1]$  and $\phi _n^y(\zeta )$  is an approximation of $y(\zeta )$  using orthonormal Bernoulli polynomials. Then the error bound of approximation can be obtained as:
\begin{equation}  \label{eq. 36}
e(y) = {\left\| {y(\zeta ) - \phi _n^y(\zeta )} \right\|_\infty } \le \frac{1}{{n!}}M
\end{equation}
where, $M = \mathop {Max}\limits_{\zeta  \in [0,1]} \phi _n^y(\zeta )y(\zeta )$.
\subsection*{Proof} See Tohidi and Kiliçman \cite{Tohidi2013} or Mahmoud \cite{Mahmoud2018}.
\subsection*{Example 1:} The Volterra integral equation 
\begin{equation} \label{eq. 37}
y(\zeta ) = 6\zeta  + 3{\zeta ^2} - \int\limits_0^\zeta  {y(\eta )d\eta } 
\end{equation}

From these theorems, it is clear that the error may be minimized to required level by including $\phi_n$  of higher degree. Furthermore, it is also obvious that the error vanishes faster with the inclusion of higher degree  $\phi_n$.

\section{Numerical Examples }
In order to discuss and establish the accuracy and effectiveness of the present method, following examples have been taken. 

\subsection*{Example 1} The Volterra integral equation 
\begin{equation}  \label{eq. 38 : Ex 1}
y(\zeta ) = 6\zeta  + 3{\zeta ^2} - \int\limits_0^\zeta  {y(\eta )d\eta } 
\end{equation}
has exact solution $ y(\zeta)=6\zeta$.

Comparing eq. (\ref{eq. 38 : Ex 1}) to standard eq. (\ref{eq. 28 : Volterra Equation}) and taking $m=5$, equations (\ref{eq. 29 : Transform of y}-\ref{eq. 32 : C in terms of Operational Matrix}) yield 
\begin{equation} \label{eq. 39: f for ex.1}
{F^T} = \left[ {4, - \frac{{3\sqrt 3 }}{2},\frac{1}{{2\sqrt 5 }},0,0,0} \right]
\end{equation}
\begin{equation} \label{eq. 40 : c for ex.1}
{C^T} = \left[ {3,\, - \sqrt 3 ,0,0,0,0} \right]
\end{equation}

Substituting eqs. (\ref{eq. 39: f for ex.1}-\ref{eq. 40 : c for ex.1}) and $\phi (\zeta ) = {[{\phi _0},\,\,{\phi _1},\,\,{\phi _2},\,...,{\phi _5}]^T}$ back into eq. (\ref{eq. 29 : Transform of y}), the exact solution $y(\zeta ) = 6\zeta $  of eq. (\ref{eq. 38 : Ex 1}) is obtained.

\subsection*{Example 2} Let us consider the Volterra integral equation of second kind
\begin{equation}  \label{eq. 41 : Ex 2}
y(\zeta ) = 1 + \zeta  - {\zeta ^2} + \int\limits_0^\zeta  {y(\eta )d\eta ,\,\,\,0}  < \zeta  < 1
\end{equation}
which has exact solution $ y(\zeta ) = 1 + 2\zeta $.

Applying the present method to eq. (\ref{eq. 41 : Ex 2}) for $m=5$  as in $example-1$, we get:
\begin{equation} \label{eq. 42: f for ex.2}
{F^T} = \left[ {\frac{7}{6},0, - \frac{1}{{6\sqrt 5 }},0,0,0} \right]
\end{equation}
\begin{equation} \label{eq. 43 : c for ex.2}
{C^T} = \left[ {2, - \frac{{393379}}{{398959\sqrt 3 }} - \frac{{1860\sqrt 3 }}{{398959}},0,0,0,0} \right]
\end{equation}

Substituting the values of $C^T$ and $F^T$ from eqs. (\ref{eq. 42: f for ex.2}-\ref{eq. 43 : c for ex.2}) and  $\phi(\zeta)$ into eq. (\ref{eq. 29 : Transform of y}), the exact solution $y(\zeta ) = 1 + 2\zeta $ of eq. (\ref{eq. 41 : Ex 2}) is obtained.

\subsection*{Example 3} Consider the following convolution integral equation
\begin{equation} \label{eq. 44 : Ex 3}
y(\zeta ) = 2 - 2{e^\zeta } + \zeta  + \frac{1}{2}{\zeta ^2} - \mathop \smallint \nolimits_0^\zeta  \left( {\zeta  - x} \right)y(x)dx
\end{equation}
having exact solution $y(\zeta ) = 1 - {e^\zeta }$.
Application of present method for $m=9$  to eq. (\ref{eq. 44 : Ex 3}),  $C^T$ and $F^T$ are obtained as :

\begin{equation} \label{eq. 45 : f for ex.3}
{F^T} = \left[ \begin{array}{l}
- \frac{{3130455131}}{{4066070400}},\frac{{2002713497}}{{2129846400\sqrt 3 }}, - \frac{{1425989}}{{7260840\sqrt 5 }}, - \frac{{578590253}}{{20766002400\sqrt 7 }},\\
- \frac{{77072}}{{116475975}}, - \frac{{1454399}}{{13214728800\sqrt {11} }}, - \frac{{445943}}{{89453548800\sqrt {13} }},\\
- \frac{{35531}}{{188278421760\sqrt {15} }},\frac{{47}}{{8132140800\sqrt {17} }}, - \frac{1}{{8821612800\sqrt {19} }}
\end{array} \right]
\end{equation}
\begin{equation} \label{eq. 46 : c for ex.3}
{C^T} = \left[ \begin{array}{l}
- 0.7182286,0.4878996, - 0.0624901, - 0.0063109,\\
- 0.0003189, - 0.101188 \times {10^{ - 4}}, - 2.177076 \times {10^{ - 7}},\\
4.442934 \times {10^{ - 9}},8.4884008 \times {10^{ - 10}}, - 7.7385133 \times {10^{ - 12}}
\end{array} \right]
\end{equation}
With help of eqs. (\ref{eq. 45 : f for ex.3}-\ref{eq. 46 : c for ex.3}),an approximate solution to eq. (\ref{eq. 44 : Ex 3}) is obtained as :

\begin{equation} \label{eq. 47 : soln to ex.3}
\begin{array}{l}
y(\zeta) = 0.002879-1.033944\zeta-0.416883\zeta^2-0.2173745\zeta^3 \\
\,\,\, -0.048615\zeta^4-0.005751\zeta^5-0.001212\zeta^6 +0.000225\zeta^7 \\
\,\,\,\,\,\, -0.000038\zeta^8-0.000002\zeta^9
\end{array}
\end{equation}

\begin{figure}[h] \label{fig.1 : plot of ex.3}
	\centering
	\includegraphics[width=0.75\linewidth]{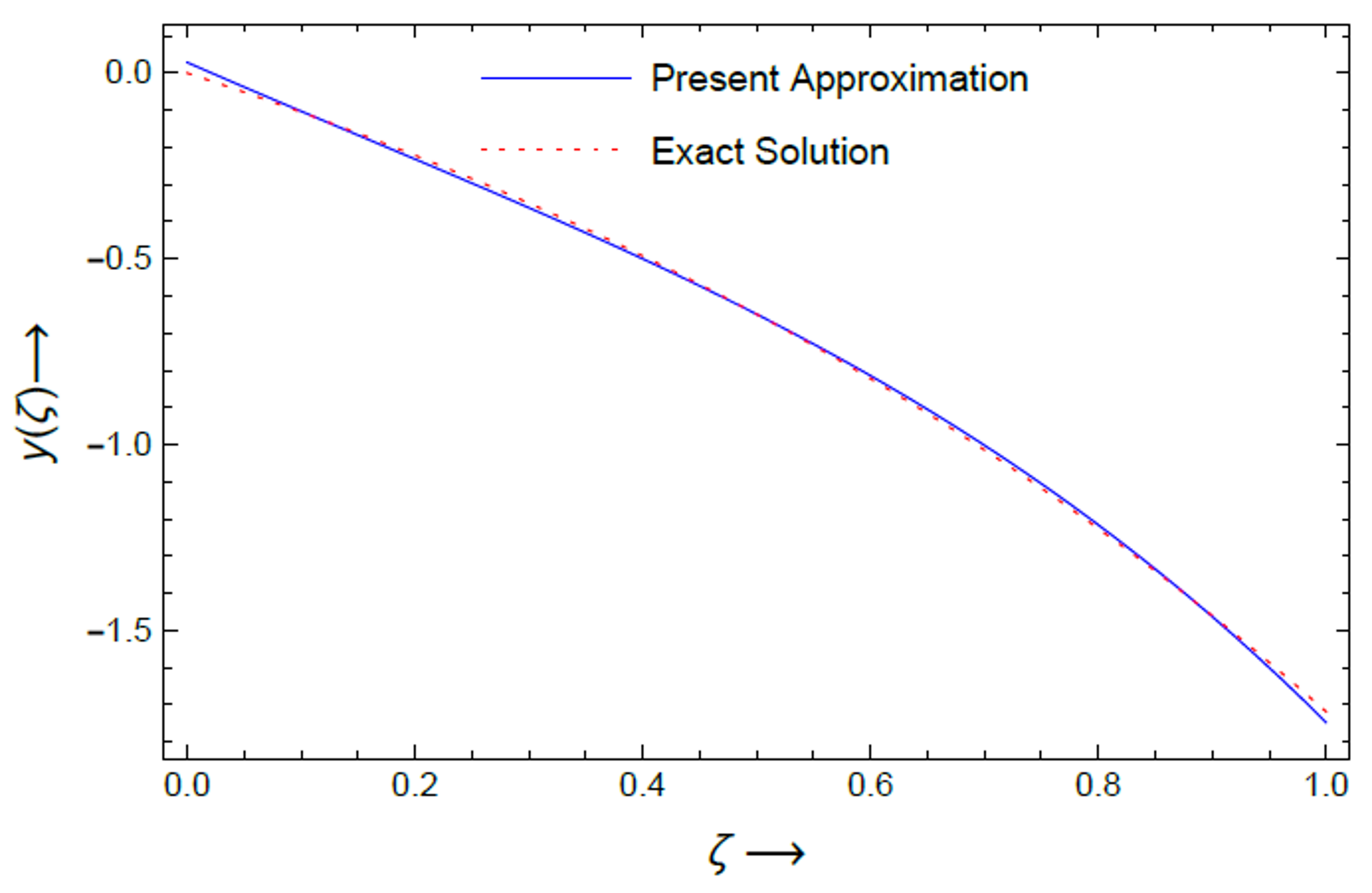}
	\caption{Comparison of Exact Solution and Approximate Solution of Example 3}
\end{figure}
\begin{figure}[h] \label{fig.2 : error of ex.3}
	\centering
	\includegraphics[width=0.7\linewidth]{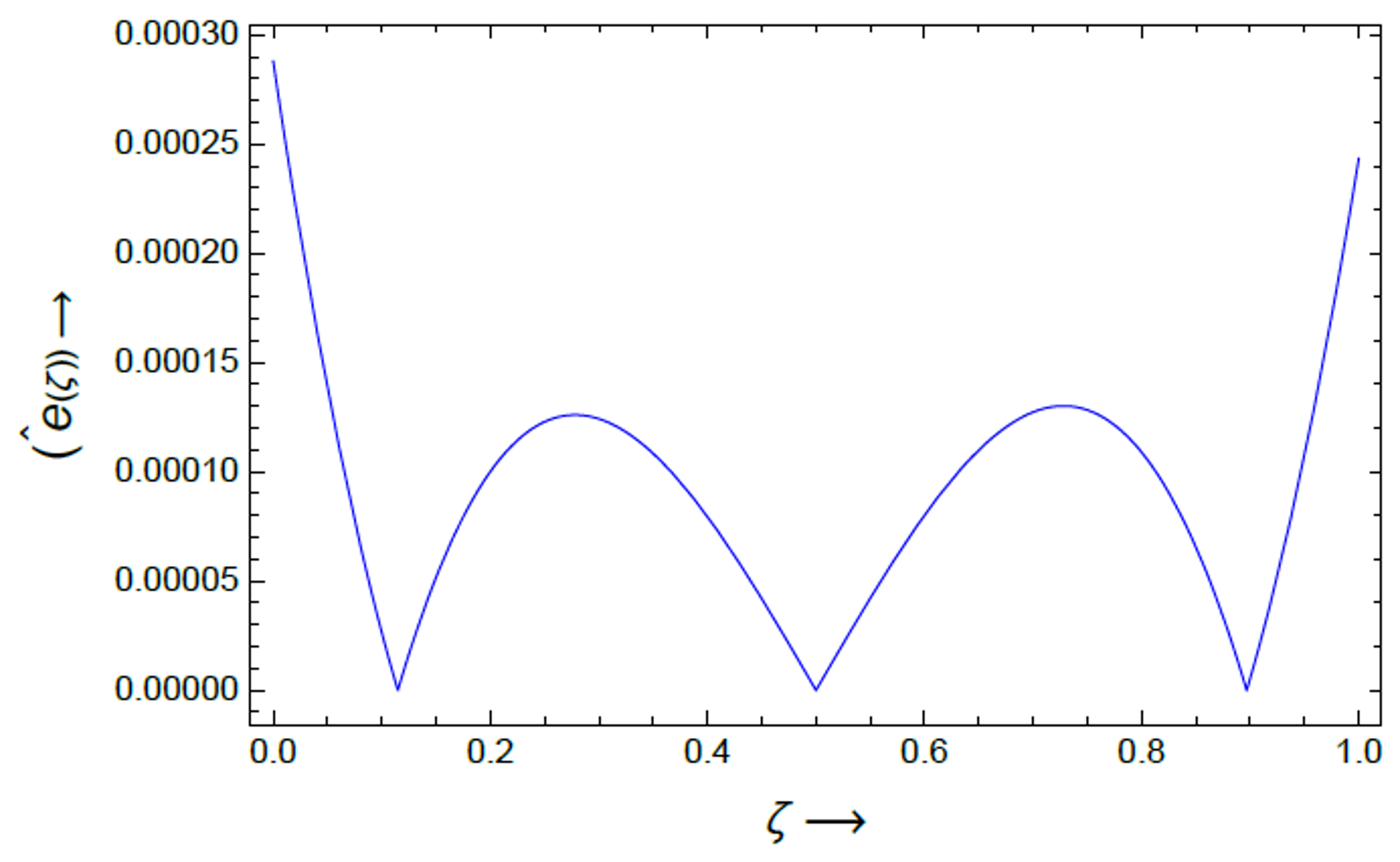}
	\caption{Absolute error, $\hat e(\zeta ),$ between exact and approximate solution of example 3}
\end{figure}

\subsection*{Example 4} Consider the following integral equation
\begin{equation} \label{eq. 48 : Ex. 4}
\begin{array}{l} 
y(\zeta ) =  - 1 - {\zeta ^2} - \frac{{{\zeta ^3}}}{3} + 2\cosh \zeta  - \sinh \zeta \\
\hspace{3.2em} + \int\limits_0^\zeta  {{{(\zeta  - \eta )}^2}y(\eta )d\eta \,\, ; \hspace{2em}\left( {0 < \zeta  < 1} \right)} 
\end{array}
\end{equation}

The exact solution of this equation is $y(\zeta ) = 1 - \sinh \zeta $.

Applying the present method for $m=9$ , we get, 
\begin{equation} \label{eq. 49 : f for Ex.4}
{F^T} = \left[ \begin{array}{l}
\frac{{1417609}}{{3628800}}, - \frac{{1025707}}{{1478400\sqrt 3 }}, - \frac{{205349}}{{1995840\sqrt 5 }}, - \frac{{322261}}{{18532800\sqrt 7 }},\\
\frac{{19}}{{54600}}, - \frac{{17}}{{5896800\sqrt {11} }},\frac{{19}}{{7257600\sqrt {13} }}, - \frac{1}{{130690560\sqrt {15} }},\\
\frac{1}{{345945600\sqrt {17} }}, - \frac{1}{{17643225600\sqrt {19} }}
\end{array} \right]
\end{equation}
\begin{equation} \label{eq. 50 : c for Ex.4}
{C^T} = \left[ \begin{array}{l}
0.45687367,\,\, - 0.33369513,\,\, - 0.0197673,\,\, - 0.00360093,\\
- 0.00010456,\,\, - 0.00001135,\,\, - 2.18846562 \times {10^{ - 7}},\,\,\\
- 1.69010986 \times {10^{ - 8}},\,\, - 2.33670044 \times {10^{ - 10}},\,\,\,0
\end{array} \right]
\end{equation}
and the solution $y(\zeta)$ is obtained as-
\begin{equation} \label{eq. 52 : Soln. of ex.4}
y(\zeta ) = {C^T}.\phi (\zeta )
\end{equation}

\begin{figure}[h] \label{fig.3 : plot of ex.4}
	\centering
	\includegraphics[width=0.75\linewidth]{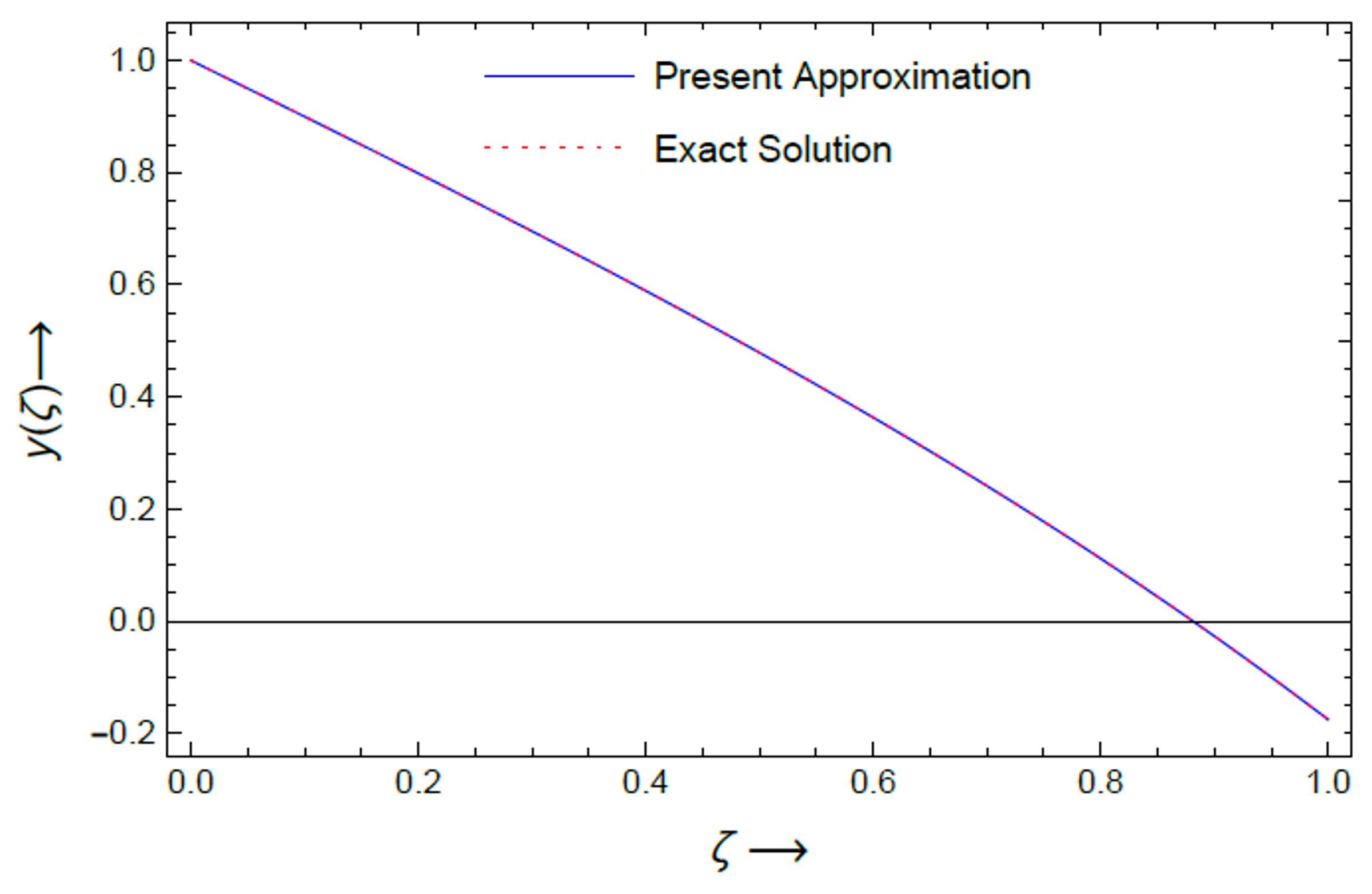}
	\caption{Comparison of Exact and Approximate Solutions of Example 4}
\end{figure}
\begin{figure}[h] \label{fig.4 : error of ex.4}
	\centering
	\includegraphics[width=0.7\linewidth]{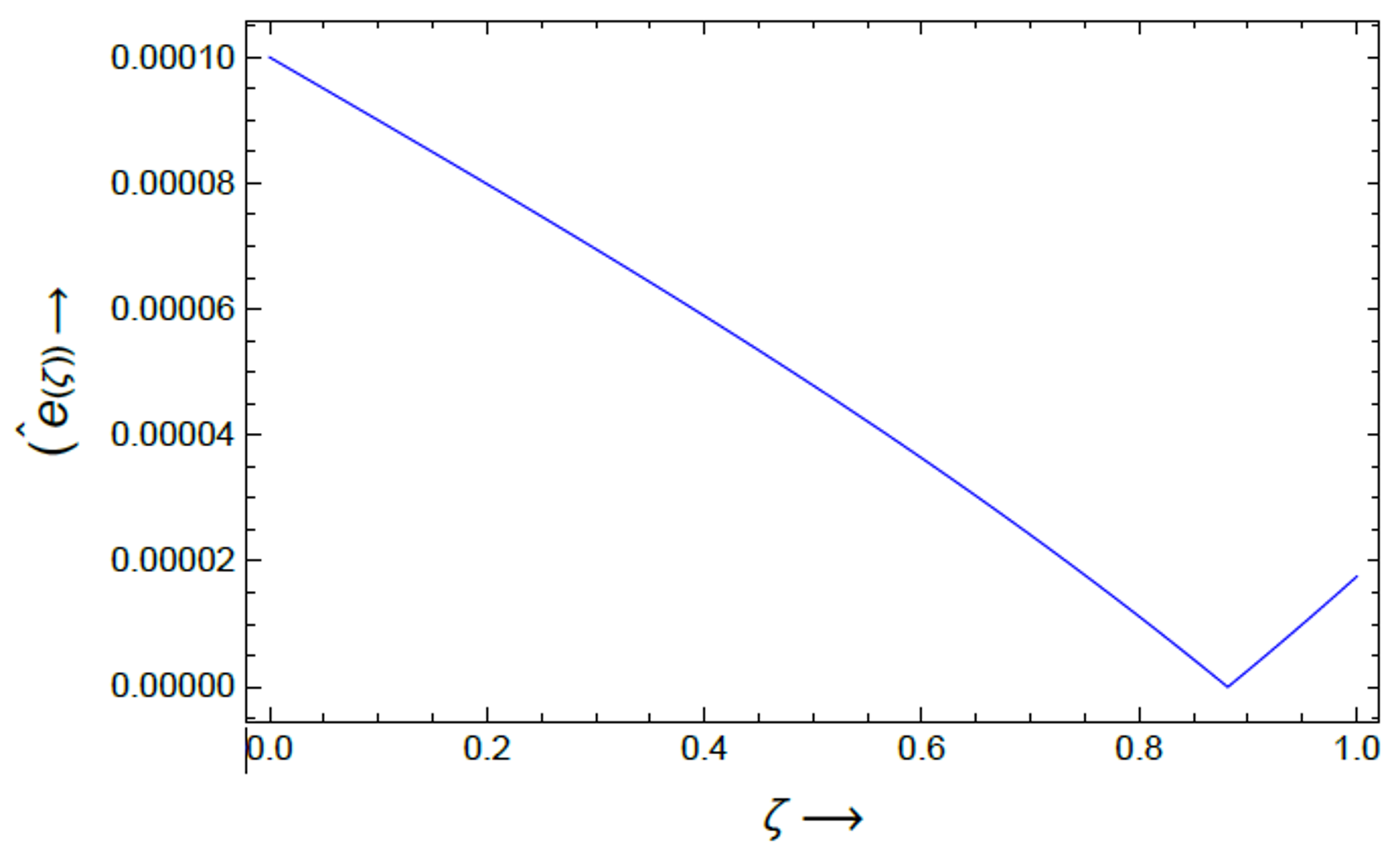}
	\caption{Absolute error, $\hat e(\zeta ),$ between exact and approximate solutions of example 4}
\end{figure}

\section{Conclusion}
In this work, we have discussed a newly developed method to find approximate solution of linear Volterra integral equations of second kind by use of Bernoulli polynomials. The process includes the derivation of an operation matrix and orthonormal polynomials. With the present process, an integral equation is converted into a system of algebraic equations with unknown coefficients, which are easily obtained with the help of coefficients generated from known part of the integral equation and operational matrix. With the help of four examples, it has been demonstrated that this method can gives either exact solution of an integral equation or an approximation in series form. Required accuracy of solution can be attained with approximation series by taking Bernoulli Polynomials of appropriate order.  

In examples 1 and 2, present method gives the exact solution with just 5 orthonormal polynomials. While, in examples 3 and 4, an approximate solution was derived with help of first nine orthonormal polynomials. The errors in examples 3 and 4 are very small in magnitude, which establish the efficacy of the present method.

The beauty of this method lies in that the method is easy for computer programming due to trigonal operational matrix, which enables to employ required number of orthonormal Bernoulli polynomials to increase the accuracy of numerical solution.

\newpage

\bibliographystyle{unsrt}

\bibliography{manuscript}

\end{document}